\documentclass[11pt,a4paper]{article}
\usepackage{epsfig,amsmath,amssymb,enumerate,mathrsfs,theorem}

\newtheorem{thm}{Theorem}[section]
\newtheorem{cor}[thm]{Corollary}
\newtheorem{lem}[thm]{Lemma}
\newtheorem{prop}[thm]{Proposition}
\newtheorem{defn}[thm]{Definition}
\newtheorem{rem}[thm]{\bf{Remark}}

\numberwithin{equation}{section}
\def\pn{\par\noindent}

\begin{document}

\title{Points at Rational Distance from the Vertices of a Unit Polygon}
\author{R. Barbara}
\date{December 15, 2009}

\maketitle

\begin{abstract}
\noindent In this paper, we investigate the existence of a point
in the plane of a unit polygon, that is at rational distance from
each vertex of the polygon. A negative answer is obtained in
almost all cases.
\end{abstract}

\pagestyle{myheadings}

\bigskip
\bigskip


\vskip 0.4 true cm

\section{\bf Introduction}

If $T$ is a unit equilateral triangle, there are points in the
plane of $T$,  that are at rational distance from the vertices of
$T$ (any vertex will do). Further, as proved in \cite{Almering}
and \cite{Berry}, the set of such points is dense in the plane of
$T$. Concerning the unit square $S$, it is not (yet) known whether
there is a point in the plane of $S$, that is at rational distance
from the corners of $S$. Results as in \cite{Berry} suggest a {\it
negative} answer, but the problem remains open.

What about the {\it unit pentagon  $P_5$} (regular pentagon with
unit side)?  Is there a point in the plane of $P_5$ that is at
rational distance from the vertices of $P_5$?

More generally, for $n\geq 3$, let $P_n$ denote the \emph{unit
$n-$gon}  (regular $n-gon$ with unit side). Consider the following
problem:

\medskip (P1) Is there a point in the plane of $P_n$ that is at rational distance from the vertices of $P_n$?

\medskip As noted, the answer to (P1) is positive if $n=3$, and it turns out that,
for $n\geq 4$, the most difficult case is indeed the case $n=4$.
In this note, we focus on the cases $n\geq 5$ and we prove the
following:

\begin{thm} \label{thm1}\mbox{ }\\
$\bullet$ For $n=5$,\qquad the answer to (P1) is NEGATIVE.\\
$\bullet$ For $n=6$,\qquad  the answer to (P1) is POSITIVE.\\
$\bullet$ For \emph{all} $n\geq 7$,  the answer to (P1) is NEGATIVE, except \emph{perhaps} if $n\in\{8,12,24\}$.
\end{thm}

The key-tool lies in the following observation: When the  answer
to (P1) is positive for a given $n\geq 3$, then, an identity as
\[
\frac{n}{4}\cot \frac{\pi}{n}=\sqrt{r_1}\pm\sqrt{r_2}\pm\cdots\pm\sqrt{r_n}
\]
must occur, where the $r_i$ are nonnegative rational numbers. But,
such identity is \emph{impossible} for $n=5$ as well as for
\emph{all} $n\geq 7$, provided that $n\neq 8, 12, 24$.

\section{\bf Preliminaries}

We start with a simple property.

\begin{prop}\label{prop1}
Let $d,m,n$ be positive integers with $d>1$ and $n=dm$. Then,
$\mathbb{Q}(\cot \frac{\pi}{d})$ and $\mathbb{Q}(\cos
\frac{2\pi}{d})$ are subfields of $\mathbb{Q}(\cot
\frac{\pi}{n})$.
\end{prop}
\vskip 0.4 true cm

\pn{\bf Proof.} $\bullet$ Set $x=\frac{\pi}{n}$ and
$y=\frac{\pi}{d}$. Then, $y=mx$. To  see why $\mathbb{Q}(\cot
y)\subset \mathbb{Q}(\cot x)$, or equivalently, $\cot y\in
\mathbb{Q}(\cot x)$, use induction on $m\geq 1$ and the identity
$\cot(m+1)x=\frac{\cot mx\cdot \cot x-1}{\cot mx+\cot x}$.
\noindent $\bullet$ Next, set $t=\cot \frac{\pi}{d}$.    From
$\cos\frac{2\pi}{d}=\frac{t^2-1}{t^2+1}$ and $t\in \mathbb{Q}(\cot
\frac{\pi}{n})$, we get, $\cos\frac{2\pi}{d}\in \mathbb{Q}(\cot
\frac{\pi}{n})$.  Hence, $\mathbb{Q}(\cos\frac{2\pi}{d})\subset
\mathbb{Q}(\cot \frac{\pi}{n})$. \hfill $\square$

\vskip 0.4 true cm Let us call a 2-group, a group in which every element has order 1 or 2. For convenience, we introduce
the
\begin{defn} We say that a real field $F$ is ''flat'' if every subfield $E$ of $F$ satisfies
\[
\text{The Galois group } G(E:\mathbb{Q}) \text{ is a 2-group}.
\]
\end{defn}

\begin{rem}  Obviously, a subfield of a flat field  is flat.
\end{rem}

\vskip 0.4 true cm
\begin{prop}\label{prop2}
Let $r_1, r_2,\ldots,r_n$ be nonnegative rational numbers. Then,
\[
\mathbb{Q}(\sqrt{r_1}\pm\sqrt{r_2}\pm\cdots\pm\sqrt{r_n}) \text{ is a \emph{flat} field}.
\]
\end{prop}

\vskip 0.4 true cm \pn{\bf Proof.} Due to the remark above, it suffices to show that $F=\mathbb{Q}(\sqrt{r_1},
\sqrt{r_2},\ldots,\sqrt{r_n})$ is a \emph{flat} field.  As quickly seen, $F:\mathbb{Q}$ is a Galois extension (of
degree $2^\nu$). We first show that  $G = G(F : \mathbb{Q})$ is a 2-group: Let $\sigma\in G$.  Then,
$\sigma(\sqrt{r_i})\in \{\pm \sqrt{r_i}\}$,   so, $\sigma\circ\sigma (\sqrt{r_i}) = \sqrt{r_i}$. As an element
$x$ in $F$ has the form $f (\sqrt{r_1} ,\sqrt{r_2}  ,\ldots ,\sqrt{r_n})$,  where  $f \in\mathbb{Q}[X_1 , X_2
,\ldots, X_n ]$, it follows easily that $\sigma\circ\sigma(x) = x$.

Since every  2-group  is abelian,   then,  $F : \mathbb{Q}$  is an \textit{abelian} extension.  Now, let $E$ be any
subfield of $F$. Since  $F:\mathbb{Q}$  is  abelian,  then,  $E : \mathbb{Q}$  is a Galois extension and the group $G(E :
\mathbb{Q})$     is isomorphic to a quotient  of  $G(F : \mathbb{Q})$.    Since a quotient of a  2-group  is a  2-group,
we see that  $G(E : \mathbb{Q})$ is a 2-group.\hfill$\square$

\begin{lem}\label{lemma}  Let  $p$  be a prime number.  Suppose that the relation  $a^2 = p(b^2 + c^2 )$
holds for some positive rational numbers $a, b, c$.   Then,  \\
$\mathbb{Q}(\sqrt{a+b\sqrt{p}}):\mathbb{Q}$ is a cyclic extension
of degree 4.
\end{lem}

\noindent\textbf{Proof. }  $\bullet$  $a + b\sqrt{p}$   is
{\small{\it{NOT}}} a square  in $\mathbb{Q}(\sqrt{p} )$:
Otherwise, for some $x, y \in \mathbb{Q}$, $a + b\sqrt{p} = (x +
y\sqrt{p} )^2$. Hence, $x^2 + py^2  = a$   and   $2xy = b$, so,
$x^2  + p\left(\frac{b}{2x} \right)^2  = a$, so, $x^2$ is a zero
of $X^2 - aX + \frac{1}{4}p b^2 = 0$. Since
$\sqrt{a^2-pb^2}=\sqrt{pc^2} =c\sqrt{p}$, it follows that  $x^2$,
hence $x$,  is irrational, a contradiction.

$\bullet$ Set $\theta = \sqrt{a + b\sqrt{p}}$. We just proved that
$\theta\notin \mathbb{Q}(\sqrt{p})$.  As further $\theta^2\in
\mathbb{Q}(\sqrt{p})$, it follows that $\theta$ has (algebraic)
degree 2  over $\mathbb{Q}(\sqrt{p} )$ and hence that $\theta$ has
degree 4 over $\mathbb{Q}$.

The irreducible polynomial of $\theta$ over $\mathbb{Q}$  is now clearly
\[
f_0 = X^4-2aX^2+(a-pb^2).
\]
The conjugates of  $\theta$ (over $\mathbb{Q}$) are: $\pm\theta$
and $\pm \mu$, where $\mu=\sqrt{a-b\sqrt{p}}$. Note that $\sqrt{p}
= \frac{1}{b}(\theta^2- a) \in\mathbb{Q}£(\theta)$. Now,
$\theta\mu = \sqrt{a^2 - pb^2}=c\sqrt{p}\in\mathbb{Q}(\theta)$.
Hence,  $\mu=\frac{c\sqrt{p}}{\theta}\in\mathbb{Q}(\theta)$.
Therefore, $\mathbb{Q}(\theta):\mathbb{Q}$ is a Galois extension
of degree 4, and hence, its Galois group $G =
G(\mathbb{Q}(\theta):\mathbb{Q})$ has order 4. Since  $f_0$ is
irreducible over $\mathbb{Q}$,   $G$  as acting on the roots of
$f_0$  is a \emph{transitive} group. In particular, for some
$\sigma\in G$,  we have,
\[
\sigma(\theta)=\mu.
\]
Claim:  $\sigma(\sqrt{p})  =  -\sqrt{p}$.  Otherwise, we must have
$\sigma(\sqrt{p})  = \sqrt{p}$,    so, $\sigma(\theta^2) =
\sigma(a + b\sqrt{p}) = a + b\sqrt{p}=\theta^2$, so,
$\sigma(\theta ) =\pm\theta$, a contradiction. Now, $\sigma(\mu) =
\sigma(\frac{c\sqrt{p}}{\theta})
=\frac{c\sigma(\sqrt{p})}{\sigma(\theta)} =
\frac{-c\sqrt{p}}{\mu}=-\theta$. Finally, $\sigma(-\theta ) = -
\mu$ and $\sigma(-\mu ) =\theta$. Hence, the action of $\sigma$ on
the roots of $f_0$ is the 4-cycle
\[
(\theta,\mu,-\theta,-\mu).
\]

As $G$ has order 4,  we conclude that $G$ is cyclic generated by $\sigma$. \hfill$\square$

\begin{prop} \label{prop3}  Each of $\mathbb{Q}(\cot\frac{\pi}{5}) :\mathbb{Q}$  and
$\mathbb{Q}(\cot\frac{\pi}{16}) :\mathbb{Q}$ is a  cyclic
extension of degree 4.
\end{prop}

\noindent \textbf{Proof. }    $\bullet$  We have
$5\cot\frac{\pi}{5}=\sqrt{25+10\sqrt{5}}$. Apply Lemma \ref{lemma}
with $p = 5$ and $(a, b, c) = (25, 10, 5)$.

$\bullet$  We have  $\cot\frac{\pi}{16}  = 1 +\sqrt{2} +
\sqrt{4+2\sqrt{2}}$.    As an exercise, check that
$\mathbb{Q}(\cot\frac{\pi}{16}) =\mathbb{Q}(\sqrt{4+2\sqrt{2}})$.
Apply Lemma \ref{lemma} with    $p = 2$ and $(a, b, c) = (4, 2,
2)$.\hfill$\square$

\begin{prop}\label{prop4}   Let  $p \geq 7$  be a prime number.
Then, $\mathbb{Q}(\cos\frac{2\pi}{p}) : \mathbb{Q}$ is a cyclic
extension of degree $\geq 3$. Further, $\mathbb{Q}(\cos
\frac{2\pi}{9}) : \mathbb{Q}$ is a cyclic extension of degree 3.
\end{prop}

\noindent\textbf{Proof. }   $\bullet$  Set  $\Omega  =
\mathbb{Q}(e^{i\frac{2\pi}{p}})$. It  is well-known that
$\Omega:\mathbb{Q}$ is a cyclic extension of degree $p -1$. Now,
$\mathbb{Q}(\cos\frac{2\pi}{p}) :\mathbb{Q}$   as a sub-extension
of    $\Omega:\mathbb{Q}$   is a cyclic extension; and it has
degree $\frac{p-1}{2}\geq 3$.

$\bullet$ Set  $\mathbb{Q}(e^{i\frac{2\pi}{9}})$. It is well-known
that   $\Omega:\mathbb{Q}$  is an  \emph{abelian} extension of
degree $\varphi(9) = 6$. Now,  $\mathbb{Q}(\cos\frac{2\pi}{9})
:\mathbb{Q}$
 as a sub-extension of an abelian  extension  is a Galois extension, so the order of its group
must be equal to its degree, that is, to $\frac{1}{2}\varphi (9) =
3$.   Since  any group of order 3  is {\small{\it{CYCLIC}}}, the
proof is complete. \hfill$\square$

\section{\bf The relation  $\dfrac{n}{4}\cot \dfrac{\pi}{n}=\sqrt{r_1}\pm\sqrt{r_2}\pm\cdots\pm\sqrt{r_n}$}

\begin{prop}\label{prop5}  Let  $n\geq 5, n\neq 6$.  Set  $\Omega  = \mathbb{Q}(\cot\frac{\pi}{n} )$.
Suppose that $\Omega$ is a flat field. Then, $n \in\{8, 12, 24\}$.
\end{prop}

\noindent\textbf{Proof.  } $\bullet$ Suppose first that $n$ is
divisible by 5. By Proposition \ref{prop1},
$\mathbb{Q}(\cot\frac{\pi}{5})$ is a subfield of $\Omega$, and, by
Proposition \ref{prop3},   the Galois group of
$\mathbb{Q}(\cot\frac{\pi}{5} ) : \mathbb{Q}$   is a \emph{cyclic
group of order 4} (hence is not a 2-group).  Therefore, $\Omega$
is {\small{\it{NOT}}} flat.

$\bullet$ Suppose next that $n$ is divisible by a prime  $p \geq
7$.     By Proposition \ref{prop1}, $\mathbb{Q}(\cos
\frac{2\pi}{p})$ is a subfield of $\Omega$, and,  by Proposition
\ref{prop4},   the Galois group of    $\mathbb{Q} (\cos
\frac{2\pi}{p}) :\mathbb{Q}$ is a \emph{cyclic group of order
$\geq 3$} (hence is not a     2-group).   Therefore, $\Omega$  is
{\small{\it{NOT}}} flat.

$\bullet$ Suppose now that $n$ is divisible by 16. By Proposition
\ref{prop1},  $\mathbb{Q}(\cot\frac{\pi}{16})$ is a subfield of
$\Omega$, and, by Proposition \ref{prop3}, the Galois group of
$\mathbb{Q}(\cot\frac{\pi}{16}) : \mathbb{Q}$ is a \emph{cyclic
group of order 4} (hence is not a 2-group). Therefore, $\Omega$ is
{\small{\it{NOT}}} flat.

$\bullet$ Suppose finally that  $n$  is divisible by  9.   By
Proposition \ref{prop1}, $\mathbb{Q}(\cos\frac{2\pi}{9})$ is a
subfield of $\Omega$, and,   by Proposition \ref{prop4},    the
Galois group of $\mathbb{Q}(\cos\frac{2\pi}{9}) :\mathbb{Q}$ is a
\emph{cyclic group of order 3} (hence is not a 2-group).
Therefore, $\Omega$ is {\small{\it{NOT}}} flat.

In conclusion, as long as we assume $\Omega$  to be flat,  $n$
cannot have a prime factor  $\geq 5$  and $n$ cannot be divisible
neither by  $2^4$   nor by  $3^2$.    Hence,  $n$  must have the
form  $n = 2^\alpha 3^\beta$ , with $\alpha\in\{0, 1, 2, 3\}$ and
$\beta\in \{0, 1\}$.  As further $n\geq 5$     and $n\neq 6$, it
remains that $n\in \{8, 12, 24\}$.\hfill$\square$

\begin{cor}\label{corprop5} Let  $n = 5$   or   $n\geq 7$,  with  $n\neq  8, 12, 24$.    Then,  an identity as
\[
\dfrac{n}{4}\cot \dfrac{\pi}{n}=\sqrt{r_1}\pm\sqrt{r_2}\pm\cdots\pm\sqrt{r_n}
\]
where the  $r_i$   are  nonnegative rational numbers,   is  \emph{impossible}.
\end{cor}

\noindent \textbf{Proof.  } Otherwise,  we would get
$\mathbb{Q}(\sqrt{r_1}\pm\sqrt{r_2}\pm\cdots\pm\sqrt{r_n}) =
\mathbb{Q}(\frac{n}{4}\cot\frac{\pi}{n} )
=\mathbb{Q}(\cot\frac{\pi}{n})$. But, by Proposition \ref{prop2},
$\mathbb{Q}(\sqrt{r_1}\pm\sqrt{r_2}\pm\cdots\pm\sqrt{r_n})$    is
a  \emph{flat}  field, whereas by  Proposition \ref{prop5},
$\mathbb{Q}(\cot\frac{\pi}{n})$ is  {\small{\it{NOT}}} a flat
field. We obtain a contradiction. \hfill$\square$

\section{\bf Proof of Theorem \ref{thm1}}

$\bullet$ For  $n = 6$,   the answer to (P1)  is  POSITIVE:   The
centroid of the unit hexagon $P_6$ is at distance  one from each
vertex.

\noindent $\bullet$ Let  $n = 5$   or   $n\geq 7$,  with  $n\neq  8, 12, 24$. We show that the answer to (P1) is
NEGATIVE. For the purpose of contradiction,     assume the existence of a point $P$ in the plane of $P_n$, that
is at rational distance from the vertices $A_1, A_2,\ldots, A_n$  of  $P_n$,  written in cyclic order. Set
$A_{n+1} = A_1$. Introduce the  $n$ triangles  $T_i = P A_i A_{i+1}$,   $i =1,\ldots,n$ (note that up to two
triangles $T_i$ might be degenerated).    Call "positive"  a triangle $T_i$ that intersects the interior of
$P_n$,  or equivalently,  such that the intersection of $T_i$ with $P_n$ has a positive area (such triangle is
non-degenerated).        Otherwise,  call  $T_i$ "negative". Note that there are always positive triangles $T_i$
(If  $P$ is interior to $P_n$ ,  all   $T_i$  are positive). Without loss of generality, we may assume that $T_1$
is positive. Now, observe the decisive properties:

(i)  If we add the areas of all positive triangles $T_i$ and then
subtract the areas of all negative triangles $T_i$   (if any),
then,  we get \emph{precisely}  the area of  $P_n$.    In other
words, we have a relation as:
\[
\text{area} ( P_n )  =  \text{area} T_1 \pm  \text{area} T_2\pm\cdots \pm\text{area} T_n.
\]

(ii)  Since every triangle $T_i$  has rational sides,  Heron's
formula $\Delta =\sqrt{s(s-a)(s-b)(s-c)}$ for the area of a
triangle shows that the area of every triangle $T_i$  has the form
$\sqrt{r_i}$,     for some nonnegative rational number $r_i$.
(Note that $\sqrt{r_i}$,  which is at most an irrational number of
degree 2,  might be rational,  even zero if $T_i$ is degenerated).

\medskip\noindent
Combining (i) and (ii),  we get,     $\text{area}(P_n )  =\sqrt{r_1}\pm\sqrt{r_2}\pm\cdots\pm\sqrt{r_n}$.

\noindent We leave it as an exercise to check that
$\text{area}(P_n) = \frac{n}{4} \cot\frac{\pi}{n}$. Finally, we
obtain
\[
 \dfrac{n}{4}\cot \dfrac{\pi}{n}=\sqrt{r_1}\pm\sqrt{r_2}\pm\cdots\pm\sqrt{r_n},
\]
 in contradiction with Corollary \ref{corprop5}.\hfill$\square$

\bigskip\noindent\textbf{Remark. } If  $P_n$  is not constructible by ruler and compasses ($\varphi(n)$ not
a power of 2), it can be shown that the (algebraic) degree of
$\frac{n}{4}\cot\frac{\pi}{n}$   over $\mathbb{Q}$  contains an
odd factor, while the degree of
$\sqrt{r_1}\pm\sqrt{r_2}\pm\cdots\pm\sqrt{r_n}$ over $\mathbb{Q}$
is a power of 2. Thus, for such $n$, the answer to (P1) is
negative. However, this will not shorten our general proof: No
decisive information is obtained for the pentagon $P_5$,  nor for
$P_n$, $n=10,15,16,17,20,30,32, \,etc$.  We even do not know
whether the constructible $P_n$, with $n$ odd, are finite or
infinite.

\bigskip\noindent \textbf{Open Problems. }  \\(1)  Solve Problem (P1) in the case  $n = 8$
(resp. $n = 12$ or $n = 24$).\\
(2) Are there points other than the centroid of the unit hexagon
$P_6$, that are  at rational distance from the vertices of  $P_6$?

\vskip 0.4 true cm

\bigskip

\bigskip
{\footnotesize \pn{\bf Roy Barbara}\; \\Lebanese University,
Faculty of Science II.\\
Fanar Campus. P.O.Box 90656. \\
Jdeidet El Metn. Lebanon.\\
{\tt Email:
roybarbara.math@gmail.com}\\
\end{document}